\documentclass[12pt]{amsart}
\usepackage[all]{xy} 

\DeclareFontEncoding{OT2}{}{} 
\newcommand{\textcyr}[1]{%
 {\fontencoding{OT2}\fontfamily{wncyr}\fontseries{m}\fontshape{n}
 \selectfont #1}}
\newcommand{\Sha}{{\mbox{\textcyr{Sh}}}}
\newcommand{\Tha}{{\mbox{\textcyr{B}}}}

\def\ra{\rightarrow}
\def \ra{\rightarrow}
\newtheorem{theorem}{Main Theorem\!\!}

\newtheorem{lemma}{Lemma}[section]

\newtheorem{proposition}[lemma]{Proposition}
\theoremstyle{definition}

\theoremstyle{remark}
\newtheorem{remark}[lemma]{Remark}

\begin{document}

\title[Cassels-Tate dual exact sequence]{A generalization of
the Cassels-Tate dual exact sequence}


\author{Cristian D. Gonz\'alez-Avil\'es}
\address{Departamento de Matem\'aticas, Universidad Andr\'es Bello,
Chile} \email{cristiangonzalez@unab.cl}

\author{Ki-Seng Tan}
\address{Department of Mathematics, National Taiwan University, Taiwan}
\email{tan@math.ntu.tw}


\thanks{C.G.-A. is partially supported by Fondecyt grant
1061209 and Universidad Andr\'es Bello grant DI-29-05/R.\\
\indent K.-S.T. is partially supported by the National Science
Council of Taiwan, NSC91-2115-M-002-001, NSC94-2115-M-002-010.}

\maketitle

\begin{abstract} We extend the well-known Cassels-Tate dual exact
sequence for abelian varieties $A$ over global fields $K$ in two
directions: we treat the $p$-primary component in the function field
case, where $p$ is the characteristic of $K$, and we dispense with
the hypothesis that the Tate-Shafarevich group of $A$ is finite.
\end{abstract}

\section{Introduction}
Let $K$ be a global field and let $m$ be a positive integer {\it
which is prime to the characteristic of $K$} (in the function field
case). Let $A$ be an abelian variety over $K$. Then there exists an
exact sequence of discrete groups
$$
0\ra\!\!\!\Sha(A)(m)\ra\!
H^{1}(K,A)(m)\ra\!\displaystyle\bigoplus_{\text{all
$v$}}H^{1}(K_{v},A)(m)\!\ra\!\!\Tha(A)(m)\ra 0,
$$
where $K_{v}$ is the henselization of $K$ at $v$,  $M(m)$ denotes
the $m$-primary component of a torsion abelian group $M$, and
$\!\Tha(A)$ is defined to be the cokernel of the localization map
$H^{1}(K,A)\ra\bigoplus_{\text{all $v$}}H^{1}(K_{v},A)$. The
Pontrjagyn dual of the preceding exact sequence is an exact sequence
of compact groups
$$
0\leftarrow\!\!\!\!\!\Sha(A)(m)^{\,*}\!\!\leftarrow
H^{1}(K,A)(m)^{\,*}\!\!\leftarrow\displaystyle\prod_{\text{all
$v$}}H^{0}(K_{v},A^{t})\widehat{\,\phantom{.}}\!\!\leftarrow\!\!
\Tha(A)(m)^{\,*}\leftarrow 0,
$$
where $A^{t}$ is the abelian variety dual to $A$ and, for any
abelian group $M$, $M\,\widehat{\,\phantom{.}}$ denotes the $m$-adic
completion $\varprojlim_{\,n} M/m^{n}$ of $M$. Now, if $\Sha(A)(m)$
is {\it finite} (or, more generally, if $\Sha(A)(m)$ contains no
nontrivial elements which are divisible by $m^{n}$ for every $n\geq
1$), then $\Sha(A)(m)^{\,*}$ and $\Tha(A)(m)^{\,*}$ are canonically
isomorphic to $\Sha(A^{t})(m)$ and
$A^{t}(K)\widehat{\,\phantom{.}}$, respectively, and the preceding
exact sequence induces an exact sequence
$$
0\leftarrow\Sha(A^{t})(m)\leftarrow
H^{1}(K,A)(m)^{\,*}\leftarrow\displaystyle\prod_{\text{all
$v$}}H^{0}(K_{v},A^{t})\widehat{\,\phantom{.}}\leftarrow
A^{t}(K)\widehat{\,\phantom{.}}\leftarrow 0
$$
which is known as {\it the Cassels-Tate dual exact sequence}
\cite{3,11}. See [9, Theorem II.5.6(b), p.247]. The aim of this
paper is to extend the isomorphism $\Tha(A)(m)^{\,*}\simeq
A^{t}(K)\widehat{\,\phantom{.}}\,$ recalled above to the case where
$m$ is divisible by the characteristic of $K$ (in the function field
case) and no hypotheses are made on $\!\!\Sha(A)$. The following is
the main result of the paper. Let $m$ and $n$ be {\it arbitrary}
positive integers. Set
$$
\text{Sel}(A^{t})_{m^{n}}=\text{Ker}\left[H^{1}(K,A^{t}_{m^{n}})\rightarrow
\bigoplus_{\text{all $v$}}H^{1}(K_{v},A^{t})\right]
$$
and
$$
T_{m}\,\textrm{Sel}(A^{t})=\varprojlim_{n}\text{Sel}(A^{t})_{m^{n}}.
$$
Then the following holds\footnote{ To see why the exact sequence of
the theorem extends the Cassels-Tate dual exact sequence recalled
above, see exact sequence (6) below and note that
$T_{m}\!\Sha(A)=T_{m}(\!\Sha(A)_{\,\text{\!$m$-div}})$ vanishes if
$\Sha(A)_{\,\text{\!$m$-div}}=0$.}:
\begin{theorem}For any positive integer $m$, there exists a natural
exact sequence of compact groups
$$
0\leftarrow\!\!\!\Sha(A)(m)^{\,*}\!\leftarrow\!
H^{1}(K,A)(m)^{\,*}\!\!\leftarrow\displaystyle\prod_{{\rm{all}}\,\,v}
H^{0}(K_{v},A^{t})\widehat{\,\phantom{.}}\leftarrow\!
T_{m}{\rm{Sel}}(A^{t})\leftarrow 0.
$$
\end{theorem}

It should be noted that a similar statement holds true if above the
henselizations of $K$ are replaced by its completions. See [9,
Remark I.3.10, p.58].

This paper grew out of questions posed to the authors by B.Poonen,
in connection with the forthcoming paper \cite{10}. We expect that
the above theorem will be useful in [op.cit.].

\section*{Acknowledgements}
K.-S.T. thanks J.Milne for the helpful outline \cite{8} which led to
a complete proof of the $p$-primary part of the Main Theorem of this
paper under the assumption that $\Sha(A)$ is finite. C.G.-A. thanks
B.Poonen for suggesting Proposition 3.3 below and, more generally,
for suggesting that the above finiteness assumption can be dispensed
with in the relevant proofs of \cite{9} if
$A(K)\widehat{\,\phantom{.}}\,$ is replaced with
$T_{m}\,\text{Sel}(A)$ throughout\footnote{ After this paper was
completed, we learned that the existence of a natural duality
between $\Tha(A)(m)$ and $T_{m}\,{\rm{Sel}}(A^{t})$ had already been
observed by J.W.S.Cassels in the case of elliptic curves over number
fields. See [3, p.153]. Therefore, the Main Theorem of this paper
may be regarded as a natural generalization of Cassels' result.}.

\section{Settings and notations}

Let $K$ be a global field and let $A$ be an abelian variety over
$K$. In the function field case, we let $p$ denote the
characteristic of $K$. All cohomology groups below are either Galois
cohomology groups or flat cohomology groups. For any non-archimedean
prime $v$ of $K$, $K_{v}$ will denote the field of fractions of the
henselization of the ring of $v$-integers of $K$. If $v$ is an
archimedean prime, $K_{v}$ will denote the completion of $K$ at $v$,
and we will write $H^{0}(K_{v},A)$ for the quotient of $A(K_{v})$ by
its identity component. Note that, for any prime $v$ of $K$, the
group $H^{1}(K_{v},A)$ is canonically isomorphic to
$H^{1}(\widehat{K}_{v},A)$, where $\widehat{K}_{v}$ denotes the
completion of $K$ at $v$. See [9, Remark I.3.10(ii), p.58]. Now let
$X$ denote either the spectrum of the ring of integers of $K$ (in
the number field case) or the unique smooth complete curve over the
field of constants of $K$ with function field $K$ (in the function
field case). In what follows, $U$ denotes a nonempty open subset of
$X$ such that $A$ has good reduction over $U$. When $N$ is a
quasi-finite flat group scheme on $U$, we endow $H^{r}(U,N)$ with
the discrete topology. Now let $m$ and $n$ be arbitrary positive
integers, and let $M$ be an abelian topological group. We will write
$M/m^{n}$ for $M/m^{n}M=M\otimes_{\mathbb{Z}}\,\mathbb{Z}/m^{n}$ and
$M\widehat{\,\,\phantom{.}}$ for the $m$-adic completion
$\varprojlim_{\,n} M/m^{n}$ of $M$. Further, we set
$\mathbb{Z}_{m}=\prod_{\,\ell\,\mid m}\mathbb{Z}_{\ell}$,
$\mathbb{Q}_{m}=\mathbb{Z}_{m}\otimes_{\mathbb{Z}}\mathbb{Q}$ and
define
$M^{\,*}=\textrm{Hom}_{\,\textrm{cts}}(M,\mathbb{Q}_{m}/\mathbb{Z}_{m})$.
Finally, the $m$-primary component of a torsion group $M$ will be
denoted by $M(m)$.

\section{Proof of the Main Theorem}
Both $A$ and its dual variety $A^{t}$ extend to abelian schemes
$\mathcal A$ and ${\mathcal A}^{t}$ over $U$ (see [2, Ch.1, \S
1.4.3]). By [5, VIII.7.1(b)], the canonical Poincar\'e biextension
of $(A^{t},A)$ by ${\Bbb G}_{m}$ extends to a biextension over $U$
of $({\mathcal A}^{t},{\mathcal A})$ by ${\Bbb G}_{m}$. Further, by
[op.cit., VII.3.6.5], (the isomorphism class of) this biextension
corresponds to a map ${\mathcal
A}^{t}\otimes^{\mathbf{L}}\!{\mathcal A}\ra {\Bbb G}_{m}[1]$ in the
derived category of the category of smooth sheaves on $U$. This map
in turn induces (see [9, p.283]) a canonical pairing
$H^{1}(U,{\mathcal A}^{t})\times H^{1}_{\textrm{c}}(U,{\mathcal
A})\ra H^{3}_{\textrm{c}}(U,{\Bbb G}_{m})\simeq \Bbb Q/ \Bbb Z$,
where the $H^{r}_{\textrm{c}}(U,{\mathcal A})$ are the cohomology
groups with compact support of the sheaf ${\mathcal A}$ defined in
[9, p.271].

\begin{remark} The smoothness of ${\mathcal A}$ implies that the
groups $H^{r}(U,{\mathcal A})$ and $H^{r}_{\textrm{c}}(U,{\mathcal
A})$ agree with the analogous groups defined for the \'etale
topology. See [9, Proposition III.0.4(d), p.272].
\end{remark}

For any positive integer $m$ and any $n\geq 1$, the above pairing
induces a pairing
\begin{equation}
H^{1}(U,{\mathcal A}^{t}_{m^{n}})\times
H^{1}_{\textrm{c}}(U,{\mathcal A})/m^{n}\ra \Bbb Q/ \Bbb Z.
\end{equation}
On the other hand, the map ${\mathcal
A}^{t}\otimes^{\mathbf{L}}\!{\mathcal A}\ra {\Bbb G}_{m}[1]$
canonically defines a map ${\mathcal A}^{t}_{m^{n}}\times{\mathcal
A}_{m^{n}}\to\mathbb{G}_{m}$, which induces a pairing
\begin{equation}
H^{1}(U,{\mathcal A}^{t}_{m^{n}})\times
H^{2}_{\textrm{c}}(U,{\mathcal A}_{m^{n}})\ra \Bbb Q/ \Bbb Z.
\end{equation}
The preceding pairing induces an isomorphism
\begin{equation}
H^{2}_{\textrm{c}}(U,{\mathcal
A}_{m^{n}})\overset{\sim}{\longrightarrow}H^{1}(U,{\mathcal
A}^{t}_{m^{n}})^{\,*}.
\end{equation}
See [9, Corollary II.3.3, p.217] for the case where $m$ prime to
$p$, and [op.cit., Theorem III.8.2, p.361] for the case where $m$ is
divisible by $p$. The pairings (1) and (2) are compatible, in the
sense that the following diagram commutes:
\begin{equation}
\xymatrix{H^{1}(U,{\mathcal A}^{t}_{m^{n}})\times
H^{1}_{\textrm{c}}(U,{\mathcal
A})/m^{n}\ar[d]^{\textrm{id}\,\times\,\partial}\ar[r]&\Bbb Q/
\Bbb Z\ar@{=}[d]\\
H^{1}(U,{\mathcal A}^{t}_{m^{n}})\times
H^{2}_{\textrm{c}}(U,{\mathcal A}_{m^{n}})\ar[r]& \Bbb Q/\Bbb Z,\\
}
\end{equation}
where $\partial\colon H^{1}_{\textrm{c}}(U,{\mathcal
A})/m^{n}\hookrightarrow H^{2}_{\textrm{c}}(U,{\mathcal A}_{m^{n}})$
is induced by the connecting homomorphism
$H^{1}_{\textrm{c}}(U,{\mathcal A})\ra
H^{2}_{\textrm{c}}(U,{\mathcal A}_{m^{n}})$ coming from the exact
sequence
$$
0\ra{\mathcal A}_{m^{n}}\ra{\mathcal
A}\overset{m^{n}}{\rightarrow}{\mathcal A}\ra 0.
$$

Now define
$$
\text{Sel}(A^{t})_{m^{n}}=\text{Ker}\left[H^{1}(K,A^{t}_{m^{n}})\rightarrow
\bigoplus_{\text{all $v$}}H^{1}(K_{v},A^{t})\right]
$$
and
$$
T_{m}\,\text{Sel}(A^{t})=\varprojlim_{n}\text{Sel}(A^{t})_{m^{n}}.
$$
By the proof of [9, Proposition I.6.4, p.92], there exists an exact
sequence
\begin{equation}
0 \to A^{t}(K)/m^{n}\to \text{Sel}(A^{t})_{m^{n}} \to
\Sha(A^{t})_{m^{n}} \to 0.
\end{equation}
Taking inverse limits, we obtain an exact sequence
\begin{equation}
0 \to A^{t}(K)\widehat{\,\phantom{.}} \to T_{m}\,\text{Sel}(A^{t})
\to T_{m}\!\Sha(A^{t}) \to 0.
\end{equation}
See [1, Proposition 10.2, p.104]. Now define\footnote{In these
definitions, the products extend over all primes of $K$, including
the archimedean primes, not in $U$.}
$$
D^{1}(U,{\mathcal A}^{t}_{m^{n}})={\rm{Ker}}\left[H^{1}(U,{\mathcal
A}^{t}_{m^{n}})\to\prod_{v\notin U}H^{1}(K_{v},A^{t})\right]
$$
and
$$
\begin{array}{rcl}
D^{1}(U,{\mathcal
A}^{t})&=&{\rm{Im}}\left[\,H^{1}_{\rm{c}}(U,{\mathcal
A}^{t})\to H^{1}(U,{\mathcal A}^{t})\,\right]\\
&=&{\rm{Ker}}\bigg[\,H^{1}(U,{\mathcal A}^{t})\to\displaystyle
\prod_{v\notin U}H^1(K_v,A^{t})\,\bigg].\\
\end{array}
$$
Note that the pairing $H^{1}(U,{\mathcal A}^{t})\times
H^{1}_{\textrm{c}}(U,{\mathcal A})\ra\Bbb Q/ \Bbb Z$ induces a
pairing $D^{1}(U,{\mathcal A}^{t})\times D^{1}(U,{\mathcal
A})\ra\Bbb Q/\Bbb Z$.

By [9, Proposition III.0.4(a), p.271] and the right-exactness of the
tensor product functor, there exists a natural exact sequence
\begin{equation}
\bigoplus_{v\notin U}H^{0}(K_{v},A)/m^{n}\to H^{1}_{\rm{c}}
(U,{\mathcal A})/m^{n}\to D^{1}(U,{\mathcal A})/m^{n}\to 0.
\end{equation}

\begin{lemma} The map $H^{1}(U,{\mathcal
A}^{t}_{m^{n}})\hookrightarrow H^{1}(K,A^{t}_{m^{n}})$ induces an
isomorphism
$$
D^{1}(U,{\mathcal A}^{t}_{m^{n}})\simeq{\rm{Sel}}(A^{t})_{m^{n}}.
$$
\end{lemma}
\begin{proof} By [9, Lemma II.5.5, p.246] and Remark 3.1 above,
the map $H^{1}(U,{\mathcal A}^{t})\!\hookrightarrow\!\!
H^{1}(K,A^{t})$ induces an isomorphism
$$
D^{1}(U,{\mathcal A}^{t})_{m^{n}}\simeq\Sha( A^{t})_{m^{n}}.
$$
Now $H^{1}(U,{\mathcal A}^{t}_{m^{n}})\to\prod_{v\in
U}H^{1}(K_{v},A^{t})$ factors through $H^{1}(U,{\mathcal
A}^{t})\to\prod_{v\in U}H^{1}(K_{v},A^{t})$, which is the zero map
(see [9, (5.5.1), p.247] and Remark 3.1 above). Consequently,
$H^{1}(U,{\mathcal A}^{t}_{m^{n}})\hookrightarrow
H^{1}(K,A^{t}_{m^{n}})$ maps $D^{1}(U,{\mathcal A}^{t}_{m^{n}})$
into ${\rm{Sel}}(A^{t})_{m^{n}}$. To prove surjectivity, we consider
the commutative diagram
\[
\xymatrix{ 0\ar[r]&{\mathcal A}^{t}(U)/m^{n}{\mathcal
A}^{t}(U)\ar[d]\ar[r] & H^{1}(U,{\mathcal
A}^{t}_{m^{n}})\ar[r]\ar[d] & H^{1}(U,{\mathcal
A}^{t})_{m^{n}}\ar[r]\ar[d]&0\\
0 \ar[r]& A^{t}(K)/m^{n}A^{t}(K)\ar[r] & H^{1}(K,A^{t}_{m^{n}})
\ar[r] &
 H^{1}(K,A^{t})_{m^{n}}\ar[r] & 0. \\
}
\]
Note that the properness of ${\mathcal A}^{t}$ over $U$ implies that
the left-hand vertical map in the above diagram is an isomorphism
(see [op.cit., p.242]). Now let $c\in{\rm{Sel}}(A^{t})_{m^{n}}$,
write $c^{\,\prime}$ for its image in $\Sha(A^{t})_{m^{n}}$ under
the map in (5) and let $\xi^{\prime}\in D^{1}(U,{\mathcal
A}^{t})_{m^{n}}\subset H^{1}(U,{\mathcal A}^{t})_{m^{n}}$ be the
pullback of $c^{\,\prime}$ under the isomorphism $D^{1}(U,{\mathcal
A}^{t})_{m^{n}}\simeq\Sha(A^{t})_{m^{n}}$ recalled above. Then the
fact that the left-hand vertical map in the above diagram is an
isomorphism implies that $\xi^{\prime}$ can be pulled back to a
class $\xi\in H^{1}(U,{\mathcal A}^{t}_{m^{n}})$ which maps down to
$c$. Clearly $\xi\in D^{1}(U,{\mathcal A}^{t}_{m^{n}})$, and this
completes the proof.\end{proof}

The following proposition generalizes [9, Theorem II.5.2(c), p.244].

\begin{proposition} There exists a canonical isomorphism
$$
(\,T_{m}\,{\rm
Sel}(A^{t}))^{\,*}\overset{\sim}{\longrightarrow}H^{2}_{\rm{c}}(U,{\mathcal
A})(m).
$$
\end{proposition}
\begin{proof} There exists a commutative diagram
\[
\xymatrix{0\ar[r] & H^{1}_{\rm{c}}(U,{\mathcal
A})/m^{n}\ar@{->}[rd]^{c}\ar[r]& H^{2}_{\rm{c}}(U,{\mathcal
A}_{m^{n}})\ar[d]^{\simeq} \ar[r]&
H^{2}_{\rm{c}}(U,{\mathcal A})_{m^{n}}\ar[r]&0\\
 &  & H^{1}(U,{\mathcal
A}^{t}_{m^{n}})^{\,*},
\\
}
\]
where the vertical map is the isomorphism (3). Clearly, the above
diagram induces an isomorphism ${\rm Coker}\,c\simeq
H^{2}_{\rm{c}}(U,{\mathcal A})_{m^{n}}$. On the other hand, there
exists a natural exact commutative\footnote{The commutativity of
this diagram follows from that of diagram (4).} diagram
\[
\xymatrix{\displaystyle\bigoplus_{v\notin U}H^{0}(K_{v},A)/m^{n}
\ar[d]^{\simeq}\ar[r]&H^{1}_{\rm{c}}(U,{\mathcal
A})/m^{n}\ar[d]^{c}\ar[r]&D^{1}(U,{\mathcal A})/m^{n}\ar[d]^{\psi}
\ar[r]&0\\
\displaystyle\bigoplus_{v\notin U}
H^{1}(K_{v},A^{t})_{m^{n}}^{\,*}\ar[r]& H^{1}(U,{\mathcal
A}^{t}_{m^{n}})^{\,*}\ar[r] &D^{1}(U,{\mathcal
A}^{t}_{m^{n}})^{\,*}\ar[r]&0,
\\
}
\]
where the top row is (7), the right-hand vertical map $\psi$ is the
composite of the natural map $D^{1}(U,{\mathcal A})/m^{n}\ra
D^{1}(U,{\mathcal A}^{t})_{m^{n}}^{\,*}$ induced by the pairing
$D^{1}(U,{\mathcal A}^{t})\times D^{1}(U,{\mathcal A})\ra\Bbb Q/\Bbb
Z$ and the natural map $D^{1}(U,{\mathcal A}^{t})_{m^{n}}^{\,*}\ra
D^{1}(U,{\mathcal A}^{t}_{m^{n}})^{\,*}$, and the left-hand vertical
map is induced by the canonical Poincar\'e biextensions of
$(A^{t},A)$ by $\mathbb{G}_{m}$ over $K_{v}$ for each $v\notin U$.
That the latter map is an isomorphism follows from [9, Remarks I.3.5
and I.3.7, pp.53 and 56, and Theorem III.7.8, p.354] and the fact
that the pairings defined in [loc.cit.] are compatible with the
pairing induced by the canonical Poincar\'e biextension (see [4,
Appendix]). The above diagram and the identification ${\rm
Coker}\,c= H^{2}_{\rm{c}}(U,{\mathcal A})_{m^{n}}$ yield an exact
sequence
$$
D^{1}(U,{\mathcal A})/m^{n}\ra D^{1}(U,{\mathcal
A}^{t}_{m^{n}})^{\,*}\ra H^{2}_{\rm{c}}(U,{\mathcal A})_{m^{n}}\ra 0
$$
Taking direct limits, we obtain an exact sequence
$$
D^{1}(U,{\mathcal A})\otimes\mathbb{Q}_{m}/\mathbb{Z}_{m}\ra
\left(\,\varprojlim D^{1}(U,{\mathcal
A}^{t}_{m^{n}})\right)^{\!*}\ra H^{2}_{\rm{c}}(U,{\mathcal A})(m)\ra
0
$$
But $D^{1}(U,{\mathcal A})\otimes\mathbb{Q}_{m}/\mathbb{Z}_{m}=0$
since $D^{1}(U,{\mathcal A})$ is torsion and
$\mathbb{Q}_{m}/\mathbb{Z}_{m}$ is divisible. Now lemma 3.2
completes the proof.
\end{proof}

By Remark 3.1 and [9, proof of Lemma II.5.5, p.247, and Proposition
II.2.3, p. 203], there exist exact sequences\footnote{In the second
exact sequence, ``$v\in U$" is shorthand for ``$v$ is a closed point
of $U$".}
$$
H^{1}(U,{\mathcal A})\overset{c_{_{
U}}}\longrightarrow\bigoplus_{v\notin U}H^{1}(K_{v},A)\ra
H^{2}_{\rm{c}}(U,{\mathcal A})
$$
and
$$
0\ra H^{1}(U,{\mathcal A})\overset{i_{_{U}}}\longrightarrow
H^{1}(K,A)\overset{\lambda_{U}}\longrightarrow\bigoplus_{v\in
U}H^{1}(K_{v},A),
$$
where $c_{_{ U}}$ and $\lambda_{U}$ are natural localization maps
and $i_{_{U}}$ is induced by the inclusion
$\textrm{Spec}K\hookrightarrow U$. If $U\subset V$ is an inclusion
of nonempty open subsets of $X$, then there exists a natural
commutative diagram
\[
\xymatrix{H^{1}(V,{\mathcal A})\ar[d]\ar[r]^{c_{_{V}}}
&\displaystyle\bigoplus_{v\notin
V}H^{1}(K_{v},A)\ar[d]\\
H^{1}(U,{\mathcal
A})\ar[r]^{c_{_{U}}}&\displaystyle\bigoplus_{v\notin
U}H^{1}(K_{v},A).
\\
}
\]
Define
$$
\Tha(A)_{U}=\textrm{coker}\!\left[c_{_{U}}\colon H^{1}(U,{\mathcal
A})\ra\bigoplus_{v\notin U}H^{1}(K_{v},A)\right],
$$
which we regard as a subgroup of $H^{2}_{\rm{c}}(U,{\mathcal A})$.
The preceding diagram shows that an inclusion $U\subset V$ of
nonempty open subsets of $X$ induces a map $\Tha(A)_{V}\ra
\Tha(A)_{U}$. Define
$$
\Tha(A)=\varinjlim
\Tha(A)_{U}=\textrm{coker}\!\left[H^{1}(K,A)\ra\bigoplus_{\text{all
$v$}}H^{1}(K_{v},A)\right],
$$
where the limit is taken over the directed family of all nonempty
open subsets $U$ of $X$ such that $A$ has good reduction over $U$,
ordered by $V\leq U$ if and only if $U\subset V$. For each $U$ as
above and every $n\geq 1$, there exists an exact sequence
$$
\bigoplus_{v\notin U}H^{1}(K_{v},A)_{m^{n}}\ra
(\!\!\Tha(A)_{U})_{m^{n}}\ra(\textrm{Im}\,c_{_U})/m^{n}.
$$
Since $\textrm{Im}\,c_{_U}$ is torsion, we conclude that there
exists a surjection
\begin{equation}
\xymatrix{\displaystyle\bigoplus_{v\notin
U}H^{1}(K_{v},A)(m)\ar@{->>}[r]&\Tha(A)_{U}(m)
\\
}
\end{equation}
On the other hand, by the proof of [9, Corollary I.6.23(b), p.111],
there exists a natural injection $T_{m}\,{\rm
Sel}(A^{t})\hookrightarrow\prod_{\,\text{all
$v$}}H^{0}(K_{v},A^{t})\widehat{\,\phantom{.}}$ and hence a
surjection
$$
\bigoplus_{\,\text{all
$v$}}(H^{0}(K_{v},A^{t})\widehat{\,\phantom{.}}\,)^{\,*}\ra(T_{m}\,{\rm
Sel}(A^{t}))^{\,*}
$$
Further, as noted in the proof of Proposition 3.3, the canonical
Poincar\'e biextensions induce an isomorphism
$$
\bigoplus_{\,\text{all
$v$}}(H^{0}(K_{v},A^{t})\widehat{\,\phantom{.}}\,)^{\,*}\simeq\,
\bigoplus_{\,\text{all $v$}}H^{1}(K_{v},A)(m),
$$
whence there exists a surjection
\begin{equation}
\xymatrix{\displaystyle\bigoplus_{\,\text{all
$v$}}H^{1}(K_{v},A)(m)\ar@{->>}[r]&(T_{m}\,{\rm Sel}(A^{t}))^{\,*}.
\\
}
\end{equation}
The maps (8) and (9) fit into a commutative diagram
\[
\xymatrix{\displaystyle\bigoplus_{\text{all $v$}}H^{1}(K_{v},A)(m)
\ar@{->>}[r]^{\textrm{\,\,\,\,\,\,(9)}}&(T_{m}\,{\rm
Sel}(A^{t}))^{\,*}\ar[r]^{\sim}&H^{2}_{\rm{c}}(U,{\mathcal
A})(m)\\
\displaystyle\bigoplus_{v\notin U}
H^{1}(K_{v},A)(m)\ar@{^{(}->}[u]\ar@{->>}[r]^{\textrm{\,\,\,\,\,\,(8)}}&
\Tha(A)_{U}(m),\ar@{^{(}-->}[u]\ar@{^{(}->}[ur]
\\
}
\]
where the isomorphism on the top row exists by Proposition 3.3.
Taking the direct limit over $U$ in the above diagram, we conclude
that there exists an isomorphism
$$
\Tha(A)(m)\overset{\sim}{\longrightarrow}(T_{m}\,{\rm
Sel}(A^{t}))^{\,*},
$$
as desired.

\begin{remark}
Recently [7, Theorem 1.2], the Cassels-Tate dual exact sequence has
been extended to 1-motives $M$ over number fields under the
assumption that the Tate-Shafarevich group of $M$ is finite. Now,
using [6, Remark 5.10], it should not be difficult to extend this
result to global function fields, provided the $p$-primary
components of the groups involved are ignored, where $p$ denotes the
characteristic of $K$. In this paper we have removed the latter
restriction when $M$ is an abelian variety, but the problem remains
for general 1-motives $M$.

\end{remark}


\begin{thebibliography}{22}
\bibitem[1]{1} Atiyah, M. and MacDonald, I.:\emph{ Introduction to
Commutative Algebra.} Addison-Wesley, Reading, MA., 1969.

\bibitem[2]{2} Bosch, S., L\"{u}tkebohmert,
W. and Raynaud, M.:\emph{ N\'eron Models.} Springer Verlag, Berlin
1989.

\bibitem[3]{3} Cassels, J.W.S:\emph{ Arithmetic on curves of genus 1.
VII. The dual
exact sequence.} J. Reine Angew. Math. {\bf{216}}, no. 1, pp.
150-158 (1964).


\bibitem[4]{4} Gonz\'alez-Avil\'es, C.D.:\emph{ Brauer groups and
Tate-Shafarevich groups.} J. Math. Sciences, Univ. Tokyo {\bf{10}},
no. 2 pp. 391-419 (2003).

\bibitem[5]{5} Grothendieck, A.:\emph{ Groupes de Monodromie en
G\'eom\'etrie Alg\'ebrique I}, S\'eminaire de G\'eom\'etrie
Alg\'ebrique du Bois Marie 1967-69 (SGA 7 I). Lecture Notes in
Math., vol. {\bf{288}}, Springer, Heidelberg, 1972.

\bibitem[6]{6} Harari, D. and Szamuely, T.:\emph{ Arithmetic duality theorems
for 1-motives.} J. reine angew. Math. {\bf{578}}, pp. 93-128 (2005).

\bibitem[7]{7} Harari, D. and Szamuely, T.:\emph{ On the arithmetic of 1-motives.}
In preparation. Available from http://www.renyi.hu/$\sim$szamuely.


\bibitem[8]{8} Milne, J.S.:\emph{ Letter to K.-S.Tan,} March 19th,
1991.

\bibitem[9]{9} Milne, J.S.:\emph{ Arithmetic Duality Theorems.}
Persp. in Math., vol. 1. Academic Press Inc., Orlando 1986.

\bibitem[10]{10} Poonen, B. and Voloch, F.:\emph{ The Brauer-Manin obstruction
for subvarieties of abelian varieties over function fields.} In
preparation.

\bibitem[11]{11} Tate, J.:\emph{ Duality theorems in Galois cohomology over
number fields.} Proceedings of the International Congress of
Mathematicians, Stockholm, 1962, pp. 288-295.


\end{thebibliography}
\end{document}